\documentclass{amsart}
\usepackage{amssymb,amsmath}

\newtheorem{theorem}{Theorem}
\newtheorem{lemma}{Lemma}

\newcommand{\bt}{\begin{theorem}}
\newcommand{\et}{\end{theorem}}
\newcommand{\bl}{\begin{lemma}}
\newcommand{\el}{\end{lemma}}
\newcommand{\beal}{\begin{align*}}
\newcommand{\enal}{\end{align*}}
\newcommand{\bq}{\begin{eqnarray*}}
\newcommand{\eq}{\end{eqnarray*}}
\newcommand{\be}{\begin{eqnarray}}
\newcommand{\ee}{\end{eqnarray}}
\newcommand{\beq}{\begin{equation}}
\newcommand{\eeq}{\end{equation}}
\newcommand{\benum}{\begin{enumerate}}
\newcommand{\eenum}{\end{enumerate}}
\newcommand{\ba}{\begin{array}}
\newcommand{\ea}{\end{array}}
\newcommand{\Z}{\ensuremath{\mathbf Z}}
\newcommand{\N}{\ensuremath{ \mathbf N }}
\newcommand{\Q}{\ensuremath{ \mathbf Q }}

\newcommand{\FF}{\ensuremath{\mathcal F}}
\newcommand{\GG}{\ensuremath{\mathcal G}}
\newcommand{\HH}{\ensuremath{\mathcal H}}
\newcommand{\supp}{\text{supp}}
\newcommand{\card}{\text{card}}

\newcommand{\pol}{$\mathcal{F} = \{f_n(q)\}_{n=1}^{\infty}$}
\newcommand{\polg}{$\mathcal{G} = \{g_n(q)\}_{n=1}^{\infty}$}
\newcommand{\polh}{$\mathcal{H} = \{h_n(q)\}_{n=1}^{\infty}$}
\newcommand{\tfe}{ the functional equation~(\ref{qc:fe})}
\newcommand{\fe}{~(\ref{qc:fe})}

\title{Quantum integers and cyclotomy}
\author{Alexander Borisov}
\address{Department of Mathematics\\ 
Pennsylvania State University\\
University Park, PA 16802}
\email{borisov@math.psu.edu}
\author{Melvyn B. Nathanson}
\address{Department of Mathematics\\
Lehman College (CUNY)\\
Bronx, New York 10468}
\email{nathansn@alpha.lehman.cuny.edu}
\thanks{The work of M.B.N. was supported
in part by grants from the NSA Mathematical Sciences Program
and the PSC-CUNY Research Award Program.}  
\author{Yang Wang}
\address{Department of Mathematics\\
Georgia Institute of Technology\\
Atlanta, GA 30332}
\email{wang@math.gatech.edu}

\thanks{Nathanson and Wang began their collaboration 
at a conference on Combinatorial and Number Theoretic Methods in Harmonic Analysis
at the Schrodinger Institute in Vienna in February, 2003.}

\begin{document}

\begin{abstract}
A sequence of functions \pol\ satisfies the functional equation 
for multiplication of quantum integers if $f_{mn}(q) = f_m(q)f_n(q^m)$
for all positive integers $m$ and $n$.
This paper describes the structure of all sequences 
of rational functions with coefficients in \Q\ 
that satisfy this functional equation.
\end{abstract}

\maketitle

\section{The functional equation for multiplication of quantum integers}
Let $\N = \{1,2,3,\ldots\}$ denote the positive integers.
For every $n \in \N$, we define the polynomial
\[
[n]_q = 1 + q + q^2 + \cdots + q^{n-1}.
\]
This polynomial is called the {\it quantum integer} $n$.
The sequence of polynomials $\{[n]_q\}_{n=1}^{\infty}$
satisfies the following functional equation:
\beq       \label{qc:fe}
f_{mn}(q) = f_m(q)f_n(q^m)
\eeq
for all positive integers $m$ and $n$.
Nathanson~\cite{nath03b} asked for a classification of all sequences
\pol\ of polynomials and of rational functions that satisfy \tfe.

The following statements are simple consequences of the functional equation.
Proofs can be found in Nathanson~\cite{nath03b}.

Let \pol\ be any sequence of functions that satisfies~(\ref{qc:fe}).
Then $f_1(q) = f_1(q)^2 = 0$ or 1.  If $f_1(q) = 0,$ then $f_n(q) = f_1(q)f_n(q) = 0$ 
for all $n \in \N$, and $\mathcal{F}$ is a trivial solution of~(\ref{qc:fe}).
In this paper we consider only nontrivial solutions of the functional equation, that is,
sequences \pol\ with $f_1(q) = 1.$

Let $P$ be a set of prime numbers, and let $S(P)$ be the multiplicative
semigroup of \N\ generated by $P$.  Then $S(P)$ consists of all integers 
that can be represented as a product of powers of prime numbers belonging to $P$.
Let \pol\ be a nontrivial solution of~(\ref{qc:fe}).
We define the support 
\[
\supp(\mathcal{F}) = \{n \in \N : f_n(q) \neq 0\}.
\]
There exists a unique set $P$ of prime numbers such that 
$\supp(\FF) = S(P)$.  Moreover, the sequence \FF\
is completely determined by the set $\{f_p(q):p\in P\}$.
Conversely, if $P$ is any set of prime numbers, and if
$\{h_p(q):p\in P\}$ is a set of functions such that
\beq          \label{qc:commute}
h_{p_1}(q)h_{p_2}(q^{p_1}) = h_{p_2}(q)h_{p_1}(q^{p_2})
\eeq
for all $p_1, p_2 \in P,$ then there exists a unique solution
\pol\ of \tfe\ such that $\supp(\FF) = S(P)$ and
$f_p(q) = h_p(q)$ for all $p \in P.$

For example, for the set $P = \{2,5,7\}$, the reciprocal polynomials
\bq
h_2(q) & = & 1 - q + q^2 \\ 
h_5(q) & = & 1 - q + q^3 - q^4 + q^5 - q^7 + q^8 \\
h_7(q) & = & 1 - q + q^3 - q^4 + q^6 - q^8 + q^9 - q^{11} + q^{12}.
\eq 
satisfy the commutativity condition~(\ref{qc:commute}). 
Since
\[
h_p(q) = \frac{[p]_{q^3}}{[p]_q} \qquad\text{for $p \in P$,}
\]
it follows that  
\beq              \label{qc:example257}
f_n(q) = \frac{[n]_{q^3}}{[n]_q} \qquad\text{for all $n \in S(P)$.}
\eeq
Moreover, $f_n(q)$ is a polynomial of degree $2(n-1)$ for all $n \in S(P).$

Let \pol\ be a solution of\tfe\ with $\supp(\FF) = S(P)$
If $P = \emptyset$, then $\supp(\FF) = \{1\}$.
It follows that $f_1(q) = 1$ and $f_n(q) = 0$ for all $n \geq 2.$
Also, for any prime $p$ and any function $h(q)$, there
is a unique solution of\tfe\ with $\supp(\FF) = S(\{ p\})$ and $f_p(q) = h(q).$
Thus, we only need to investigate solutions of~(\ref{qc:fe}) for
$\card(P) \geq 2.$

If \pol\ and \polg\ are solutions of\fe\ with $\supp(\FF) = \supp(\GG),$
then, for any integers $d$, $e$, $r$, and $s$, 
the sequence of functions \polh, where
\[
h_n(q) = f_n(q^r)^dg_n(q^s)^e,
\]
is also a solution of\tfe\ with $\supp(\HH) = \supp(\FF).$
In particular, if \pol\ is a solution of\fe, then \polh\ is another solution of\fe,
where 
\[
h_n(q) = \left\{
\ba{ll}
1/f_n(q)& \text{if $n \in \supp(\FF)$}\\
0 & \text{if $n \not\in \supp(\FF)$.}
\ea
\right.
\]
The functional equation also implies that 
\beq          \label{qc:mn}
f_m(q)f_n(q^m) = f_n(q)f_m(q^n)
\eeq
for all positive integers $m$ and $n$, and
\beq            \label{qc:prod}
f_{m^k}(q) = \prod_{i=0}^{k-1}f_m(q^{m^i}).
\eeq

Let \pol\ be a solution in rational functions of\tfe\
with $\supp(\FF) = S(P).$
Then there exist a completely multiplicative arithmetic function
$\lambda(n)$ with support $S(P)$ and rational numbers $t_0$ and $t_1$ 
with $t_0(n-1) \in \Z$ and $t_1(n-1) \in \Z$ for all $n \in S(P)$ such that,
for every $n \in S(P)$, 
we can write the rational function $f_n(q)$ uniquely in the form
\beq             \label{qc:ratform}
f_n(q) = \lambda(n)q^{t_0(n-1)}\frac{u_n(q)}{v_n(q)},
\eeq
where $u_n(q)$ and $v_n(q)$ are monic polynomials 
with nonzero constant terms, and
\[
\deg(u_n(q)) - \deg(v_n(q)) = t_1(n-1) \quad\text{ for all $n \in \supp(\FF)$.}
\]
For example, let $P$ be a set of prime numbers with $\card(P) \geq 2.$
Let $\lambda(n)$ be a completely multiplicative arithmetic function
with support $S(P)$, and let $t_0$ be a rational number
such that $t_0(n-1) \in \Z$ for all $n \in S(P)$.
Let $R$ be a finite set of positive integers and $\{t_r\}_{r\in R}$ a set of integers. 
We construct a sequence \pol\ of rational functions as follows:
For $n \in S(P),$ we define 
\beq           \label{qc:finalform}
f_n(q) = \lambda(n)q^{t_0(n-1)} \prod_{r\in R} [n]_{q^r}^{t_r} .
\eeq
For $n \not\in S(P)$ we set $f_n(q) = 0$.
Then $\prod_{r\in R} [n]_{q^r}^{t_r}$ is a quotient of monic polynomials 
with coefficients in \Q\ and nonzero constant terms.
The sequence \pol\ satisfies\tfe, and $\supp(\FF) = S(P)$.

We shall prove that every solution of\tfe\ in rational functions with
coefficients in \Q\ is of the form~(\ref{qc:finalform}).
This provides an affirmative answer to Problem 6 in~\cite{nath03b}
in the case of the field \Q.

\section{Roots of unity and solutions of the functional equation}
Let $K$ be an algebraically closed field, and let $K^*$ denote the multiplicative group
of nonzero elements of $K.$
Let ${\Gamma}$ denote the group of roots of unity in $K^*$,
that is,
\[
{\Gamma} = \{\zeta \in K^* : \zeta^n = 1 \text{ for some $n \in N$}\}.
\]
Since ${\Gamma}$ is the torsion subgroup of $K^*,$ 
every element in $K^*\setminus {\Gamma}$ has infinite order.
We define the {\it logarithm group}
\[
L(K) = K^*/{\Gamma},
\]
and the map
\[
L:K^* \rightarrow L(K)
\]
by
\[
L(a) = a{\Gamma} \text{ for all $x \in K^*$}.
\] 
We write the group operation in $L(K)$ additively:
\[
L(a) + L(b) = a{\Gamma} + b{\Gamma} = ab{\Gamma} = L(ab).
\]

\bl
Let $K$ be an algebraically closed field, and $L(K)$ its logarithm group.
Then $L(K)$ is a vector space over the field \Q\ of rational numbers.
\el

\begin{proof}
Let $a \in K^*$ and $m/n \in \Q.$  Since $K$ is algebraically closed, 
there is an element $b \in K^*$ such that
\[ 
b^n = a^m.
\]
We define
\[
\frac{m}{n}L(a) = L(b).
\]
Suppose $m/n = r/s \in \Q$, and that
\[
c^s = a^r
\]
for some $c \in K^*$.
Since $ms = nr,$ it follows that
\[
c^{ms} = a^{mr} = b^{nr} = b^{ms},
\]
and so $c/b \in {\Gamma}.$  Therefore,
\[
\frac{m}{n}L(a) = L(b) = b{\Gamma} = c{\Gamma} = L(c) = \frac{r}{s}L(a),
\]
and $(m/n)L(a)$ is well-defined.
It is straightforward to check that $L(K)$ is a \Q-vector space.
\end{proof}

\bl         \label{qc:lemma:log}
Let $P$ be a set of primes, $\card(P) \geq 2,$ and let
$S(P)$ be the multiplicative semigroup generated by $P$.
For every integer $m \in S(P)\setminus \{1\}$ there is an integer $n \in S(P)$
such that $\log m$ and $\log n$ are linearly independent over \Q.
Equivalently, for every integer $m \in S(P)\setminus \{1\}$ 
there is an integer $n \in S(P)$ such that 
there exist integers $r$ and $s$ with
$m^r = n^s$ if and only if $r = s = 0.$
\el

\begin{proof}
If $m = p^k$ is a prime power, let $n$ be any prime in $P \setminus \{p\}.$
If $m$ is divisible by more than one prime, let $n$ be any prime in $P$.
The result follows immediately from the Fundamental Theorem of Arithmetic. 
\end{proof}

Let $K$ be a field.
A {\it function on $K$} is a map
$f:K \rightarrow K \cup \{\infty\}$. 
For example, $f(q)$ could be a polynomial or a rational function 
with coefficients in $K.$
We call $f^{-1}(0)$ the set of {\it zeros} of $f$ and 
$f^{-1}(\infty)$ the set of {\it poles} of $f$.

\bt           \label{qc:theorem:zeropole}
Let $K$ be an algebraically closed field.
Let \pol\ be a sequence of functions on $K$ that satisfies \tfe.
Let $P$ be the set of primes such that $\supp(\FF) = S(P).$
If $\card(P) \geq 2$ and if, for every $n \in \supp(\FF)$, 
the function $f_n(q)$ has only finitely many zeros and only finitely many poles, 
then every zero and pole of $f_n(q)$ is either 0 or a root of unity.
\et

\begin{proof}
The proof is by contradiction.  
Let ${\Gamma}$ be the group of roots of unity in $K$.
Suppose that 
\[
f_n(a) = 0 \text{ for some $n \in \supp(\FF)$
and $a \in K^*\setminus {\Gamma} .$}
\]  
By Lemma~\ref{qc:lemma:log}, there is an integer $m \in S(P)$ such that
$\log m$ and $\log n$ are linearly independent over \Q.
Since $a$ has infinite order in the multiplicative group $K^*$
and $f_n^{-1}(0)$ is finite, 
there are positive integers $k$ and $M = m^k$
such that $a^{M}$ is not a zero of the function $f_n(q).$ 
By~(\ref{qc:mn}), we have
\[
f_M(q)f_n(q^M) = f_n(q)f_M(q^n).
\]
Therefore,
\[
f_M(a)f_n(a^M) = f_n(a)f_M(a^n) = 0.
\]
Since $f_n(a^M) \neq 0,$ it follows from~(\ref{qc:prod}) that
\[
0 = f_M(a) = f_{m^k}(a) = \prod_{i=0}^{k-1} f_m(a^{m^i}),
\]
and so
\[
f_m(a^{m^i}) = 0 \text{ for some $i$ such that $0 \leq i \leq k-1.$}
\]
Let
\[
b = a^{m^i}.
\]
Then 
\[
f_m(b) = 0,
\]
\[
b \in K^*\setminus {\Gamma},
\]
and
\beq        \label{qc:a}
L(b) = m^iL(a)
\eeq

Since $f_m^{-1}(0)$ is finite,
there are positive integers $\ell$ and $N = n^{\ell}$
such that $z^N$ is not a zero of $f_m(q)$ 
for every $z \in f_m^{-1}(0)$ with $z \in K^*\setminus {\Gamma}.$
Since $K$ is algebraically closed, we can choose $c \in K$ such that 
\[
c^N = b.
\]
Then 
\[
f_m(c) \neq 0,
\]
\[
c \in K^*\setminus {\Gamma},
\]
and
\beq        \label{qc:b}           
NL(c) = L(b).
\eeq
Again applying~(\ref{qc:mn}), we have
\[
f_m(q)f_N(q^m) = f_N(q)f_m(q^N) 
\]
and so
\[
f_m(c)f_N(c^m) = f_N(c)f_m(c^N) = f_N(c)f_m(b) = 0.
\]
It follows that
\[
0 = f_N(c^m) = f_{n^{\ell}}(c^m) = \prod_{j=0}^{\ell-1} f_n(c^{mn^j}),
\]
and so
\[
f_n(c^{mn^j}) = 0 \text{ for some $j$ such that $0 \leq j \leq \ell-1.$}
\]
Let 
\[
a' = c^{mn^j}.
\]
Then
\[
f_n(a') = 0,
\]
\[
a' \in K^*\setminus {\Gamma},
\]
and
\beq           \label{qc:c}
L(a') = mn^jL(c)
\eeq
Combining~(\ref{qc:a}),~(\ref{qc:b}), and~(\ref{qc:c}),
we obtain
\[
L(a') = \frac{mn^j}{N}L(b) = \frac{m^{i+1}}{n^{\ell-j}}L(a),
\] 
that is,
\beq                 \label{qc:a'}
L(a') = \frac{m^{i'}}{n^{j'}}L(a),
\text{ where $1 \leq i' \leq k \text{ and } 1 \leq j' \leq \ell.$}
\eeq

What we have accomplished is the following:  Given an element $a \in f_n^{-1}(0)$
that is neither 0 nor a root of unity, we have constructed another element 
$a'\in f_n^{-1}(0)$ that is also neither 0 nor a root of unity, and that 
satisfies~(\ref{qc:a'}).  Iterating this process,
we obtain an infinite sequence of such elements.  
However, the number of zeros of $f_n(q)$ is finite, and so the elements 
in this sequence cannot be pairwise distinct.  It follows that there is an element
\[
a  \in f_n^{-1}(0) \setminus \left({\Gamma} \cup \{0\}\right)
\]
such that
\[
L(a) = \frac{m^r}{n^s} L(a),
\]
where $r$ and $s$ are positive integers.
Then
\[
a^{n^s}{\Gamma} = L\left( a^{n^s}\right) = n^sL(a) 
= m^rL(a) = L\left( a^{m^r}\right) = a^{m^r}{\Gamma}.
\]
Since $a$ is not a root of unity,
it follows that
\[
m^r = n^s,
\]
which contradicts the linear independence of $\log m$ and $\log n$ over \Q.
Therefore, the zeros of the functions $f_n(q)$ belong to ${\Gamma} \cup \{0\}$ 
for all $n \in \supp(\FF).$

Replacing the sequence $\mathcal{F} = \{f_n(q)\}_{n\in \supp(\FF)}$
with $\mathcal{F}'= \{1/f_n(q)\}_{n\in \supp(\FF)}$,
we conclude that the poles of the functions $f_n(q)$ also belong to ${\Gamma} \cup \{0\}$ 
for all $n \in \supp(\FF).$
This completes the proof.
\end{proof}

\section{Rational solutions of the functional equation}
In this section we shall completely classify sequences 
of rational functions with rational coefficients
that satisfy the functional equation for quantum multiplication.

For $k \geq 1,$ let $\Phi_k(q)$ denote the $k$th cyclotomic polynomial.  Then
\[
F_k(q) = q^k -1 = \prod_{d|k} \Phi_d(q)
\]
and
\beq           \label{qc:cyclopol}
\Phi_k(q) = \prod_{d|k}F_d(q)^{\mu(k/d)},
\eeq
where $\mu(k)$ is the M\" obius function.
Let $\zeta$ be a primitive $d$th root of unity.  
Then $F_k(\zeta) = 0$ if and only if $d$ is a divisor of $k$.
We define 
\[
F_0(q) = \Phi_0(q) = 1.
\]
Note that
\beq         \label{qc:pol}
F_k(q) = q^k-1 = (q-1)(1+q+ \cdots + q^{k-1}) = F_1(q)[k]_q
\eeq
for all $k \geq 1.$

A {\em multiset} $U = (U_0 ,\delta)$ consists of a finite set $U_0$ 
of positive integers and a function $\delta: U_0 \rightarrow \N$.
The positive integer $\delta(u)$ is called the {\em multiplicity} of $u$.
Multisets $U = (U_0 ,\delta)$ and $U' = (U'_0 ,\delta')$ are equal if
$U_0 = U'_0$ and $\delta(u) = \delta'(u)$ for all $u \in U_0.$
Similarly, $U \subseteq U'$ if $U_0 \subseteq U'_0$ and $\delta(u) \leq \delta'(u)$
for all $u \in U_0.$
The multisets $U$ and $U'$ are {\em disjoint} if $U_0 \cap U'_0 = \emptyset.$
We define
\[
\prod_{u\in U} f_u(q) = \prod_{u\in U_0} f_u(q)^{\delta(u)}
\]
and
\[
\max(U) = \max(U_0).
\]
If $U_0 = \emptyset,$ then we set $\max(U) = 0$ and $\prod_{u\in U} f_u(q) = 1.$

\bl                                                    \label{qc:lemma:unique}
Let $U$ and $U'$ be multisets of positive integers.
Then
\beq           \label{qc:unique}
\prod_{u\in U} F_u(q) = \prod_{u'\in U'} F_{u'}(q),
\eeq
if and only if $U = U'.$
\el

\begin{proof}
Let $k = \max(U \cup U').$  Let $\zeta$ be a primitive $k$th root of unity.
If $k \in U',$ then
\[
\prod_{u\in U} F_u(\zeta) = \prod_{u'\in U'} F_{u'}(\zeta) = 0,
\]
and so $k \in U.$  
Dividing~(\ref{qc:unique}) by $F_k(q),$ reducing the multiplicity of $k$ in 
the multisets $U$ and $U'$ by 1, and continuing inductively, we obtain $U = U'.$
\end{proof}

Let \pol\ be a nontrivial solution of\tfe, where
$f_n(q)$ is a rational function with rational coefficients for all $n \in \supp(\FF)$.
Because of the standard representation~(\ref{qc:ratform}), we can assume that
\[
f_n(q) = \frac{u_n(q)}{v_n(q)},
\]
where $u_n(q)$ and $v_n(q)$ are monic polynomials 
with nonzero constant terms.
By Theorem~\ref{qc:theorem:zeropole}, the zeros of the polynomials $u_n(q)$ and $v_n(q)$
are roots of unity, and so we can write
\[
f_n(q) = \frac{\prod_{u\in U'_n}\Phi_u(q)}{\prod_{v\in V'_n}\Phi_v(q)},
\]
where $U'_n$ and $V'_n$ are disjoint multisets of positive integers.
Applying~(\ref{qc:cyclopol}), we replace each cyclotomic polynomial 
in this expression with a quotient of polynomials of the form $F_k(q)$.
Then
\beq                   \label{qc:Frep}
f_n(q) = \frac{\prod_{u\in U_n} F_u(q)}{\prod_{v\in V_n} F_u(q)},
\eeq
where $U_n$ and $V_n$ are disjoint multisets of positive integers.
Let
\[
f_n(q) = \frac{\prod_{u\in U_n} F_u(q)}{\prod_{v\in V_n} F_v(q)}
= \frac{\prod_{u'\in U'_n} F_{u'}(q)}{\prod_{v'\in V'_n} F_{v'}(q)},
\]
where $U_n$ and $V_n$ are disjoint multisets of positive integers
and $U'_n$ and $V'_n$ are disjoint multisets of positive integers.
Then
\[
\prod_{u\in U_n \cup V'_n} F_u(q) = \prod_{v\in U'_n \cup V_n} F_v(q).
\]
By Lemma~\ref{qc:lemma:unique}, we have the multiset identity
\[
U_n \cup V'_n = U'_n \cup V_n.
\]
Since $U_n \cap V_n = \emptyset,$ it follows that $U_n \subseteq U'_n$
and so $U_n = U'_n.$  Similarly, $V_n = V'_n.$
Thus, the representation~(\ref{qc:Frep}) is unique.

We introduce the following notation for the {\em dilation} of a set:  
For any integer $d$ and any set $S$ of integers,
\[
d\ast S = \{ds: s \in S\}.
\]

\bl           \label{qc:lemma:mp}
Let \pol\ be a nontrivial solution of\tfe\
with $\supp(\FF) = S(P)$, where $\card(P) \geq 2.$
Let
\[
f_n(q) = \frac{\prod_{u\in U_n} F_u(q)}{\prod_{v\in V_n} F_v(q)}
\]
and $U_n$ and $V_n$ are disjoint multisets of positive integers.
For every prime $p \in P,$ let 
\[
m_p = \max(U_p \cup V_p).
\]
There exists an integer $r$ such that $m_p = rp$ for every $p \in P.$
Moreover, either $m_p \in U_p$ for all $p \in P$ or $m_p \in V_p$ for all $p \in P$.
\el

\begin{proof}
Let $p_1$ and $p_2$ be prime numbers in $P$, and let
\[
\frac{m_{p_1}}{p_1} \geq \frac{m_{p_2}}{p_2}.
\]
Equivalently,
\[
p_2m_{p_1} \geq p_1m_{p_2}.
\]
Applying functional equation~(\ref{qc:mn}) with $m = p_1$ and $n = p_2,$
we obtain
\[
\frac{\prod_{u\in U_{p_1}} F_u(q)}{\prod_{v\in V_{p_1}} F_v(q)}
\frac{\prod_{u\in U_{p_2}} F_u(q^{p_1})}{\prod_{v\in V_{p_2}} F_v(q^{p_1})}
=\frac{\prod_{u\in U_{p_2}} F_u(q)}{\prod_{v\in V_{p_2}} F_v(q)}
\frac{\prod_{u\in U_{p_1}} F_u(q^{p_2})}{\prod_{v\in V_{p_1}} F_v(q^{p_2})},
\]
where
\[
U_{p_1} \cap V_{p_1} = U_{p_2} \cap V_{p_2} = \emptyset.
\]
The identity
\[
F_n(q^m) = \left(q^m\right)^n - 1 = q^{mn}-1 = F_{mn}(q),
\]
implies that
\begin{align*}
\frac{\prod_{u\in U_{p_1}\cup p_1\ast U_{p_2}} F_u(q)}
{\prod_{v\in V_{p_1}\cup p_1\ast V_{p_2}} F_v(q)} 
= &
\frac{\prod_{u\in U_{p_1}} F_u(q)}{\prod_{v\in V_{p_1}} F_v(q)}
\frac{\prod_{u\in p_1\ast U_{p_2}} F_u(q)}{\prod_{v\in p_1\ast V_{p_2}} F_v(q)} \\
= & \frac{\prod_{u\in U_{p_2}} F_u(q)}{\prod_{v\in V_{p_2}} F_v(q)}
\frac{\prod_{s\in p_2\ast U_{p_1}} F_u(q)}{\prod_{t\in p_2\ast V_{p_1}} F_v(q)} \\
= & \frac{\prod_{u\in U_{p_2}\cup p_2\ast U_{p_1}} F_u(q)}
{\prod_{v\in V_{p_2}\cup p_2\ast V_{p_1}} F_v(q)}.
\end{align*}
By the uniqueness of the representation~(\ref{qc:Frep}),
it follows that
\[
U_{p_1} \cup (p_1\ast U_{p_2}) \cup V_{p_2} \cup (p_2\ast V_{p_1})
= U_{p_2} \cup (p_2\ast U_{p_1}) \cup V_{p_1} \cup (p_1\ast V_{p_2}).
\]

Recall that 
\[
m_{p_1} = \max\left(U_{p_1} \cup V_{p_1}\right).
\]
If
\[
m_{p_1} \in U_{p_1},
\]
then
\[
p_2m_{p_1} \in p_2\ast U_{p_1}
\]
and so
\[
p_2m_{p_1} \in U_{p_1} \cup (p_1\ast U_{p_2}) \cup V_{p_2} \cup (p_2\ast V_{p_1}).
\]
However, 
\benum
\item[(i)]
$p_2m_{p_1} \not\in U_{p_1}$ 
since $p_2m_{p_1} > m_{p_1} = \max\left(U_{p_1} \cup V_{p_1}\right),$
\item[(ii)]
$p_2m_{p_1} \not\in p_2\ast V_{p_1}$ since 
$m_{p_1} \in U_{p_1}$ and $U_{p_1}\cap V_{p_1} = \emptyset,$ 
\item[(iii)]
$p_2m_{p_1} \not\in V_{p_2}$ 
since $p_2m_{p_1} \geq p_1m_{p_2} > m_{p_2}= \max\left(U_{p_2} \cup V_{p_2}\right)$.
\eenum
If $p_2m_{p_1} > p_1m_{p_2} = \max\left( p_1\ast U_{p_2} \right),$ then
$p_2m_{p_1} \not\in p_1\ast U_{p_2}$.
This is impossible, and so
\[
p_2m_{p_1} = p_1m_{p_2}\in p_1\ast U_{p_2},
\]
\[
m_{p_2}\in U_{p_2},
\]
and
\[
\frac{m_{p_1}}{p_1} = \frac{m_{p_2}}{p_2} = r \quad\text{ for all $p_1,p_2 \in P.$}
\]

Similarly, if $m_{p_1}\in V_{p_1}$ for some $p_1 \in P,$
then $m_{p_2}\in V_{p_2}$ for all $p_2 \in P$.
This completes the proof.
\end{proof}

\bt          \label{qc:theorem:structure}
Let \pol\ be a sequence of rational functions with coefficients in \Q\
that satisfies\tfe.  
If $\supp(\FF) = S(P),$ where $P$ is a set of prime numbers and $\card(P) \geq 2,$
then there are
\benum
\item[(i)]
a completely multiplicative arithmetic function $\lambda(n)$ with support $S(P)$, 
\item[(ii)]
a rational number $t_0$ such that $t_0(n-1)$ is an integer for all $n \in S(P)$, 
\item[(iii)]
a finite set $R$ of positive integers and a set $\{t_r\}_{r\in R}$ of integers
\eenum
such that
\beq         \label{qc:structure}
f_n(q) = \lambda(n) q^{t_0(n-1)}\prod_{r\in R} [n]_{q^r}^{t_r} \quad\text{ for all $n \in \supp(\FF)$.}
\eeq
\et

\begin{proof}
It suffices to prove~(\ref{qc:structure}) for all $p \in P.$
Recalling the representation~(\ref{qc:ratform}), 
we only need to investigate the case 
\[
f_p(q) = \frac{\prod_{u\in U_p} F_u(q)}{\prod_{v\in V_p} F_v(q)},
\]
where $U_p$ and $V_p$ are disjoint multisets of positive integers.
Let $m_p = \max(U_p\cup V_p).$
By Lemma~\ref{qc:lemma:mp}, there is a nonnegative integer $m$ such that
$m_p = mp$ for all $p \in P.$  We can assume that $m_p \in U_p$ for all $p \in P.$  

The proof is by induction on $m$.  
If $m = 0$, then $U_p = V_p = \emptyset$ and $f_p(q) = 1$ for all $p \in P,$
hence~(\ref{qc:structure}) holds with $R = \emptyset.$

Let $m = 1,$ and suppose that $m_p = p \in U_p$ for all $p \in P.$
Then 
\[
f_p(q) = \frac{\prod_{u\in U_p} F_u(q)}{\prod_{v\in V_p} F_v(q)} \\
= \frac{(q^p-1)\prod_{u\in U'_p} F_u(q)}{\prod_{v\in V_p} F_v(q)}.
\]
Since $q^p-1 = F_1(q)[p]_q,$ we have
\bq
g_p(q) & = & \frac{f_p(q)}{[p]_q} \\
& = & \frac{(q^p-1)\prod_{u\in U_p\setminus\{p\}} F_u(q)}{[p]_q\prod_{v\in V_p} F_v(q)} \\
& = & \frac{F_1(q)\prod_{u\in U_p\setminus\{p\}} F_u(q)}{\prod_{v\in V_p} F_v(q)} \\
& = & \frac{\prod_{u\in U'_p} F_u(q)}{\prod_{v\in V'_p} F_v(q)},
\eq
where $U'_p \cap V'_p = \emptyset$. 
The sequence of rational functions \polg\ is also a solution of\tfe,
and either $\max(U'_p \cup V'_p) = 0$ for all $p \in P$
or $\max(U'_p \cup V'_p) = p$ for all $p \in P.$

If $\max(U'_p \cup V'_p) = p$ for all $p \in P$,
then we construct the sequence \polh\ of rational functions
\[
h_n(q) = \frac{g_n(q)}{[n]_q} = \frac{f_n(q)}{[n]_q^2}.
\]
Continuing inductively, we obtain a positive integer $t$ such that
\[
f_n(q) = [n]_q^t   \quad\text{ for all $n \in \supp(\FF).$}
\]
Thus,~(\ref{qc:structure}) holds in the case $m = 1.$

Let $m$ be an integer such that
the Theorem holds whenever $m_p < mp$ for all $p \in P,$ 
and let \pol\ be a solution of\tfe\ with $\supp(\FF) = S(P)$
and $m_p = mp$ and $m_p \in U_p$ for all $p \in P.$  
The sequence \polg\ with
\[
g_n(q) = \frac{f_n(q)}{[n]_{q^r}}
\]
is a solution of\tfe.
Since
\[
F_{rp}(q) = q^{rp}-1 
= \left(q^r-1\right)\left( 1 + q^r + \cdots + q^{r(p-1)} \right) 
= F_r(q)[p]_{q^r},
\]
it follows that 
\bq
g_p(q) 
& = & \frac{(q^{m_p}-1)\prod_{u\in U_p\setminus\{m_p\}}F_u(q)}{[p]_{q^r}\prod_{v\in V_p} F_v(q)} \\
& = & \frac{(q^{mp}-1)\prod_{u\in U_p\setminus\{mp\}} F_u(q)}{[p]_{q^r}\prod_{v\in V_p} F_v(q)} \\
& = & \frac{F_r(q)\prod_{u\in U_p\setminus\{mp\}} F_u(q)}{\prod_{v\in V_p} F_v(q)} \\
& = & \frac{\prod_{u\in U'_p} F_u(q)}{\prod_{v\in V'_p} F_v(q)},
\eq
where $U'_p \cap V'_p = \emptyset$, and $\max(U_{p'}\cup V_{p'}) \leq mp.$

If $\max(U_{p'}\cup V_{p'}) = mp,$ then $mp \in U'_p.$
We repeat the construction with
\[
h_n(q) = \frac{g_n(q)}{[n]_{q^r}} = \frac{f_n(q)}{[n]_{q^r}^2}.
\]
Continuing this process, we eventually obtain a positive integer $t_r$ 
such that the sequence of rational functions
\[
\left\{\frac{f_n(q)}{[n]_{q^r}^{t_r}}\right\}_{n=1}^{\infty}
\]
satisfies\tfe,
and
\[
\frac{f_p(q)}{[p]_{q^r}^{t_r}}
= \frac{\prod_{u\in U'_p} F_u(q)}{\prod_{v\in V'_p} F_v(q)},
\]
where $U'_p \cap V'_p = \emptyset$ and $\max(U'_p \cup V'_p) < mp.$
It follows from the induction hypothesis there is a finite set $R$
of positive integers and a set $\{t_r\}_{r\in R}$ of integers such that
\[
f_n(q) = \prod_{r\in R} [n]_{q^r}^{t_r} \quad\text{ for all $n \in \supp(\FF)$.}
\]
This completes the proof.
\end{proof}

There remain two related open problems.
First, we would like to have a simple criterion to determine when a sequence 
of rational functions satisfying\tfe\ is actually a sequence of polynomials.
It is sufficient that all of the integers $t_r$ in the 
representation~(\ref{qc:structure}) be nonnegative,
but the example in~(\ref{qc:example257}) shows that this condition is not necessary. 

Second, we would like to have a structure theorem for rational function solutions
and polynomial solutions to\tfe\ with coefficients in an arbitrary field,
not just the field of rational numbers.

\providecommand{\bysame}{\leavevmode\hbox to3em{\hrulefill}\thinspace}
\providecommand{\MR}{\relax\ifhmode\unskip\space\fi MR }
\providecommand{\MRhref}[2]{%
  \href{http://www.ams.org/mathscinet-getitem?mr=#1}{#2}
}
\providecommand{\href}[2]{#2}


\begin{thebibliography}{1}

\bibitem{nath03b}
M.~B. Nathanson, \emph{A functional equation arising from multiplication of
  quantum integers}, J. Number Theory, to appear.

\end{thebibliography}
\end{document}